# Traffic dynamics and optimal control in a city served by ride-sourcing vehicles

Wenbo Fan [*]

School of Transportation and Logistics, Southwest Jiaotong University, Chengdu, China; [*] Corresponding email: wbfan@swjtu.edu.cn

**Abstract**

This paper presents an interactive bathtub model for describing the traffic dynamics of ride-sourcing vehicles including non-shared taxis and ride-pooling cars. A city with a network of undifferentiated streets and solely served by ride-sourcing services is assumed to facilitate the modeling, isolate the congestion contribution, and accordingly develop control strategies. The proposed model is parsimonious with only input information of the total lane length of the network, the in-flux of demand, and the travel distance distributions. The output of the model, however, captures not only the traffic dynamics of vehicles but also the dynamic states of passengers in ride-pooling services in terms of the total number, the remaining travel distances, and the queue of unmatched requests at any system time. Useful system metrics can be exploited for use of the authorities to monitor, predict, and control the traffic, as well as for the TNCs to determine the fleet sizes, dispatch vehicles, and measure the service productivity. For illustration, we propose a robust control rule to manage the traffic efficiently and avoid gridlock, and also present time-varying ride-pooling sizes to eliminate the queue of unmatched requests. Numerical examples demonstrate the effectiveness of the proposed model and the control and operation strategies.

**Keywords**: traffic dynamics; bathtub model; ride-sourcing services; gridlock; optimal control

## 1. Introduction

In the past decades, ride-sourcing services, with the aid of mobile and wireless communication technologies, have witnessed an explosive growth in many cities all over the world. Evidences of their impacts on traffic are, however, discouraging. For instance, a recent report indicated that ride-sourcing service providers (also called the Traffic Network Companies, TNCs) had been adding considerable amount of new vehicle miles traveled (VMT) in major cities of US, and the number is 976 million miles of driving on New York City streets from 2013 to 2017 (Schaller Consulting, 2018); such an increase in VMT may reach as high as 85% in the Denver area (Henao, 2017). Similar findings were also reported in other areas/countries (Hawkins, 2019; Nie, 2017). Evidently, without proper control, ride-sourcing services may worsen the traffic congestion, which would also backfire on the TNCs.

Proper control should be based on accurate descriptions of the traffic dynamics of ride-sourcing vehicles. Resorting to the literature, we unfortunately find no existing models ready for this task. The majority studies on ride-sourcing services were focused on the operational tactics, e.g., demand-supply matching, static and dynamic vehicle routing, and pricing, with a variety of service types and practical constraints (See recent overviews in Mourad et al., 2019; Agatz et al., 2012; and Furuhata et al., 2013). An exception is the work by Daganzo and Ouyang (2019), who proposed a general model describing the states of ride-sourcing vehicles in all service stages, i.e., requests-vehicle matching, picking-up/collecting, and delivering. Their model is, however, a static one, and also ignores the congestion effect in the network by assuming a constant speed. Alternatively, simulation-based methods are indeed applicable, which had been used in studying the dynamic fleet sizing of shared vehicles (Fagnant and Kockelman, 2016), the dynamic interaction with public transit (Mo et al, 2020), and the joint design with public transit (Pinto et al., 2019; Gurumurthy et al., 2020). These works are, however, very time and cost consuming due to the complex structure of many simulation modules and huge data requirement.

Most recently, a breakthrough in the so-called bathtub model provides us an excellent tool to our problem. Jin (2020) formalized Vickrey's bathtub model (VBM) (Vickrey, 1991; 2020) into a generalized form (GBM) with two major contributions: (i) Vickrey's assumption of the negative exponential distribution of trip distances was relaxed to account for any distributions; and (ii) the traffic dynamics were tracked and evolved in terms of the total number of vehicles (trips) in the network at any time and the distribution of their remaining distances to travel. The GBM and VBM recognize the network performance of any speed-density relations, e.g., the Macroscopic Fundamental Diagram (MDF) among others, and thus are capable of capturing the phenomenon of hyper-congestion and gridlock. With these capabilities, the bathtub models have a very appealing property of parsimoniousness: the only input information is the total lane length from the supply side, and the in-flux and trip distance distributions of the demand. In the seminal work, Vickrey formulated the VBM as an ordinary derivative equation (ODE) using the number of vehicles as the state variable. In retrospect, similar ODEs were independently discovered in Agnew (1976), Mahmassani and Herman (1984), and Daganzo (2007). None of them, however, explicitly clarified the three premises of VBM: (i) an undifferentiated network; (ii) a speed-density relation; and (iii) a negative exponential distribution of trip distances. Bathtub models of other forms had been continuously studied by Arnott and his co-authors (Arnott, 2013; Arnott et al., 2016; Arnott and Buli, 2018), who assumed identical travelers with the same constant trip distance, built delay-differential equations (DDEs), but shown some difficulties in solving them. Applications of VBM can be seen in Small and Chu (2003) and Fosgerau (2015) for solving departure time user equilibrium.

The above studies, however, are limited to the traffic dynamics of private cars only. Therefore, we intend to make an extension of bathtub model (i.e., GBM in particular) to the case of ride-sourcing vehicles, which become an increasingly important contributor to traffic congestion. The main contributions of this paper are three-fold: (i) an interactive bathtub model is established for describing the dynamics of ride-sourcing vehicles in three states: idling, picking-up/collecting, and delivering; (ii) an robust control strategy is proposed to manage the traffic efficiently and avoid gridlock; and (iii) for ride-pooling services, the optimal ride-pooling size is presented with the promise of eliminating queues of unmatched requests.

The rest of this paper is organized as follows. Next section begins the formulations for the dynamic equations of non-shared taxis, followed by the solution method in Section 3. Section 4 extends the models to the case of ride-pooling/sharing with participants more than one. Numerical examples are presented in Section 5 demonstrating the effectiveness of the proposed model and solution method. Last section draws main conclusions and discusses future research directions.

## 2. Models

Consider a city (generally represents a certain urban area) with a maze of undifferentiated streets, of which the total lane kilometers is $L$ (lane-km). During the study period (e.g., morning peak hours), the demand is characterized by a continuous in-flux function, $f(t)$ (trips/h). For the transportation system, we consider only ride-sourcing services (i.e., non-shared taxis and ride-pooling cars) in this study and treated other transportation modes (e.g., private cars, buses, metro) as exogenous given and fixed as if the city is solely served by the former. In doing so, we isolate the dynamics for the ride-sourcing vehicles in congested traffic and accordingly develop control strategies to avoid gridlock.

As the beginning, we start from modeling non-shared taxis; extension is then done to ride-pooling services with ride-pooling size $c > 1$. For the sake of consistence with previous studies, we adopt most notations from Jin (2020) and Daganzo and Ouyang (2019). A full list of notations is presented in Table A1 of Appendix A.

## 2.1 Interactive bathtub model for non-shared taxis

The non-shared taxis can be classified into three groups according to their service states: (i) the idle ones waiting to be matched with emerging requests; (ii) the matched ones running to pick up the assigned requests; and (iii) the occupied ones delivering passengers to their destinations. We follow the notation used in Daganzo and Ouyang (2019) to denote the numbers of the three types of taxis in the system at time $t$: $n_{00}(t), n_{01}(t)$, and $n_{10}(t)$ (vehicles); and thus the total number of taxis is $N(t) = n_{00}(t) + n_{01}(t) + n_{10}(t)$ (vehicles).

The non-idle/active vehicles are characterized by PDF functions of the desired travel distance distributions, $\tilde{\varphi}_{01}(t,x)$ and $\tilde{\varphi}_{10}(t,x)$ (km$^{-1}$), respectively. These characteristics are dictated by ride-haling requests and treated as exogenously given, which can be obtained via surveys, mining from historical data, or simulation.[1] The service states of the active vehicles are described by PDF functions of the remaining travel distance distributions, $\varphi_{01}(t,x)$ and $\varphi_{10}(t,x)$ (km$^{-1}$), respectively. These PDFs satisfy:
(i) $\tilde{\varphi}_{01}(t,x), \tilde{\varphi}_{10}(t,x), \varphi_{01}(t,x), \varphi_{10}(t,x) \geq 0$; (ii) $\int_0^\infty \tilde{\varphi}_{01}(t,x)dx = \int_0^\infty \tilde{\varphi}_{10}(t,x)dx = \int_0^\infty \varphi_{01}(t,x)dx = \int_0^\infty \varphi_{10}(t,x)dx = 1$; and (iii) $\tilde{\varphi}_{01}(t,\infty) = \tilde{\varphi}_{10}(t,\infty) = \varphi_{01}(t,\infty) = \varphi_{10}(t,\infty) = 0$.

At any time $t$, the traffic dynamics of the active vehicles satisfy the following conservation laws:

$$n_{01}(t+\Delta t)\varphi_{01}(t+\Delta t, x)\Delta x = n_{01}(t)\varphi_{01}(t, x+v(t)\Delta t)\Delta x + a_{00}(t)\Delta t \tilde{\varphi}_{01}(t,x)\Delta x \quad (1)$$

$$n_{10}(t+\Delta t)\varphi_{10}(t+\Delta t, x)\Delta x = n_{10}(t)\varphi_{10}(t, x+v(t)\Delta t)\Delta x + p_{01}(t)\Delta t \tilde{\varphi}_{10}(t,x)\Delta x \quad (2)$$

where the left-hand side (LHS) of (1,2) are the numbers of the matched/collecting vehicles and of the occupied/delivering vehicles, respectively, whose remaining travel distances are between $x$ and $x + \Delta x$ at $t + \Delta t$; the first term in the right-hand side (RHS) of (1,2) are the numbers of two types of vehicles whose remaining travel distances were between $x + v(t)\Delta t$ and $x + v(t)\Delta t + \Delta x$ at $t$, which would be shorten by $v(t)\Delta t$ and become in between $x$ and $x + \Delta x$ at $t + \Delta t$; and the second term in the RHS of (1,2) are the numbers of two types of vehicles that newly emerge during $t$ and $t + \Delta t$. The $a_{00}(t)$ (vehicles/h) is the in-flux of collecting vehicles and equivalently the out-flux of idle vehicles, and $a_{00}(t) = f(t)$ if no control/constraint is imposed on supply/matching of ride-sourcing vehicles; otherwise under a control strategy, $a_{00}(t) = \min(f(t), \bar{a})$, and $\bar{a}$ functions as a cap of the maximum ride-sourcing supply/request-vehicle-matching rate. The $p_{01}(t)$ is the in-flux of delivering vehicles, i.e., the out-flux of the collecting ones that finished picking up passengers.

---

[1] Immense data are nowadays available in the TNCs, which can be required to provide such information of travel distance distributions of e-hailing vehicles.

Define the densities of active vehicles with remaining travel distance being $x$ at $t$ as $k_{01}(t,x) \equiv n_{01}(t)\varphi_{01}(t,x)$ and $k_{01}(t,x) \equiv n_{10}(t)\varphi_{10}(t,x)$, respectively, (1) and (2) can be rewritten as follows by eliminating $\Delta x$ from all terms,

$$k_{01}(t+\Delta t, x) = k_{01}(t, x+v(t)\Delta t) + a_{00}(t)\Delta t \tilde{\varphi}_{01}(t,x) \quad (3)$$

$$k_{10}(t+\Delta t, x) = k_{10}(t, x+v(t)\Delta t) + p_{01}(t)\Delta t \tilde{\varphi}_{10}(t,x) \quad (4)$$

Replacing $k_{01}(t, x+v(t)\Delta t)$ with $k_{01}(t,x) + \frac{\partial}{\partial x}k_{01}(t,x)v(t)\Delta t$ in (3) and $k_{10}(t, x+v(t)\Delta t)$ with $k_{10}(t,x) + \frac{\partial}{\partial x}k_{10}(t,x)v(t)\Delta t$ in (4), and dividing $\Delta t$ on both sides, the above equations become:

$$\frac{\partial}{\partial t}k_{01}(t,x) - \frac{\partial}{\partial x}k_{01}(t,x)v(t) = a_{00}(t)\tilde{\varphi}_{01}(t,x) \quad (5)$$

$$\frac{\partial}{\partial t}k_{10}(t,x) - \frac{\partial}{\partial x}k_{10}(t,x)v(t) = p_{01}(t)\tilde{\varphi}_{10}(t,x) \quad (6)$$

where $p_{01}(t)$ can be derived from the in-flux and out-flux conservation: $p_{01}(t)\Delta t = \int_{x=0}^{v(t)\Delta t} k_{01}(t,x)dx \cong k_{01}(t,0)v(t)\Delta t$, which yields $p_{01}(t) \cong k_{01}(t,0)v(t)$ (vehicles/h).

The two partial derivative equations (PDEs) of (5,6) describe the traffic dynamics of active vehicles over time $t$ in terms of their densities and remaining travel distances. To solve them, we need know $v(t)$ in the system. In general, $v(t)$ can be expressed as an unimodal function of the average vehicular density $\rho(t) \equiv \frac{N(t)}{L}$: $v(t) = V(\rho(t))$, or $v(t) = V\left(\frac{N(t)}{L}\right)$, i.e.,

$$v(t) = V\left(\frac{n_{00}(t) + \int_0^\infty k_{01}(t,x)dx + \int_0^\infty k_{10}(t,x)dx}{L}\right) \quad (7)$$

where $n_{00}(t)$ can be described using the following dynamics:

$$\dot{n}_{00}(t) = s(t) + d_{10}(t) - a_{00}(t) \quad (8)$$

where $s(t)$ is the exogenous supply rate of ride-sourcing vehicles, which can be positive, zero, or negative; $d_{10}(t)$ is the in-flux from the delivering vehicles that finished the delivery mission; and thus can be estimated by $d_{10}(t) \cong k_{10}(t,0)v(t)$. Note that $s(t)$ can be set as $a_{00}(t) - d_{10}(t)$, such that $\dot{n}_{00}(t) = 0$, and there is no oversupply or deficiency of ride-sourcing fleet.[2]

Overall, the traffic dynamics of the ride-sourcing vehicles are completely described by the PDEs (5,6), the speed-density function (7), and the ODE (8), which can be summarized as following interactive bathtub model:

---

[2] Considering that some idle taxies may cruise on streets but some may be parked, a discount coefficient $\alpha \in [0,1]$ can be applied to $n_{00}(t)$ in (5). For the sake of conservation, we simply set $\alpha = 1$ here.

$$\dot{n}_{00}(t) = s(t) + k_{10}(t,0)V\left(\frac{n_{00}(t)+\int_0^\infty k_{01}(t,x)dx+\int_0^\infty k_{10}(t,x)dx}{L}\right) - a_{00}(t) \quad (9a)$$

$$\frac{\partial}{\partial t}k_{01}(t,x) - \frac{\partial}{\partial x}k_{01}(t,x)V\left(\frac{n_{00}(t)+\int_0^\infty k_{01}(t,x)dx+\int_0^\infty k_{10}(t,x)dx}{L}\right) = a_{00}(t)\tilde{\varphi}_{01}(t,x) \quad (9b)$$

$$\frac{\partial}{\partial t}k_{10}(t,x) - \frac{\partial}{\partial x}k_{10}(t,x)V\left(\frac{n_{00}(t)+\int_0^\infty k_{01}(t,x)dx+\int_0^\infty k_{10}(t,x)dx}{L}\right) = k_{01}(t,0)v(t)\tilde{\varphi}_{10}(t,x) \quad (9c)$$

where three unknown variables, $n_{00}(t)$, $k_{01}(t,x)$, $k_{10}(t,x)$, for the three types of vehicles can be solved given certain strategies of $s(t)$ and $a_{00}(t)$ and the information of $f(t)$, $\tilde{\varphi}_{01}(t,x)$, and $\tilde{\varphi}_{10}(t,x)$.

2.2 Derivations of useful system metrics

2.2.1 Number of active vehicles, $n_{01}(t)$, $n_{10}(t)$

Although the number of active vehicles at any time $t$ can be directly obtain by $n_{01}(t) = \int_0^\infty k_{01}(t,x)dx$ and $n_{10}(t) = \int_0^\infty k_{10}(t,x)dx$, they may also be derived and represented by more computation-efficient forms as follows. Define $K_{01}(t,x) \equiv \int_{y=x}^\infty k_{01}(t,y)dy$ as the number of active pickup vehicles with remaining trip distance no less than $x$; and $K_{10}(t,x) \equiv \int_{y=x}^\infty k_{10}(t,y)dy$ as the number of active delivery vehicles with remaining trip distance no less than $x$. Integrating both sides of (5,6) yields,

$$\frac{\partial}{\partial t}K_{01}(t,x) - v(t)\frac{\partial}{\partial x}K_{01}(t,x) = a_{00}(t)\tilde{\Phi}_{01}(t,x) \quad (10)$$

$$\frac{\partial}{\partial t}K_{10}(t,x) - v(t)\frac{\partial}{\partial x}K_{10}(t,x) = p_{01}(t)\tilde{\Phi}_{10}(t,x) \quad (11)$$

where $\tilde{\Phi}_{01}(t,x) \equiv \int_{y=x}^\infty \tilde{\varphi}_{01}(t,y)dy$ and $\tilde{\Phi}_{10}(t,x) \equiv \int_{y=x}^\infty \tilde{\varphi}_{10}(t,y)dy$ are the cumulative distribution functions (CDFs) of active vehicles with remaining travel distance no less than $x$.

Note that $\frac{\partial}{\partial t}K_{01}(t,x) - v(t)\frac{\partial}{\partial x}K_{01}(t,x) = \frac{d}{dt}K_{01}(t,\tilde{x}_0 - z(t))$ and $\frac{\partial}{\partial t}K_{10}(t,x) - v(t)\frac{\partial}{\partial x}K_{10}(t,x) = \frac{d}{dt}K_{10}(t,\tilde{x}_0 - z(t))$, where $\tilde{x}_0$ represents the desired travel distance of any vehicle entering the system at the initial time $t = 0$; and $z(t) \equiv \int_0^t v(s)ds$ is defined as the characteristic distance traveled (by such a vehicle) in the system since $t = 0$. The meaning of $\tilde{x}_0 - z(t)$ is the remaining travel distance $x$. Thus, we have:

$$\frac{d}{dt}K_{01}(t,x) = a_{00}(t)\tilde{\Phi}_{01}(t,\tilde{x}_0 - z(t)) \quad (12)$$

$$\frac{d}{dt}K_{10}(t,x) = p_{01}(t)\tilde{\Phi}_{10}(t,\tilde{x}_0 - z(t)) \quad (13)$$

Solving the ODEs (12,13) yields $K_{01}(t,x) = K_{01}(t,\tilde{x}_0) + \int_0^t a_{00}(s)\widetilde{\Phi}_{01}(t,\tilde{x}_0 - z(s))ds$ and $K_{10}(t,x) = K_{10}(t,\tilde{x}_0) + \int_0^t p_{01}(s)\widetilde{\Phi}_{10}(t,\tilde{x}_0 - z(s))ds$, in which replacing $\tilde{x}_0 = x + z(t)$ gives us:

$$K_{01}(t,x) = K_{01}(t, x+z(t)) + \int_0^t a_{00}(s)\widetilde{\Phi}_{01}(t, x+z(t) - z(s))ds \quad (14)$$

$$K_{10}(t,x) = K_{10}(t, x+z(t)) + \int_0^t p_{01}(s)\widetilde{\Phi}_{10}(t, x+z(t) - z(s))ds \quad (15)$$

Knowing that $K_{01}(t,0) = n_{01}(t)$ and $K_{10}(t,0) = n_{10}(t)$; thus, from (14,15) we obtain the number of active vehicles at any time $t$:

$$n_{01}(t) = K_{01}(t, z(t)) + \int_0^t a_{00}(s)\widetilde{\Phi}_{01}(t, z(t) - z(s))ds \quad (16)$$

$$n_{10}(t) = K_{10}(t, z(t)) + \int_0^t p_{01}(s)\widetilde{\Phi}_{10}(t, z(t) - z(s))ds \quad (17)$$

Note that (14-17) are directly computed by variables of time $t$ without the need of integrating $k_{01}(t,x)$ and $k_{10}(t,x)$ over $x$ for each $t$, which implies a higher computation efficiency.

### 2.2.2 Other metrics

Define $w(t)$ as the number of requests waiting for assignment to ride-sourcing vehicles, the dynamics of $w(t)$ can be expressed by:

$$\dot{w}(t) = f(t) - a_{00}(t) \quad (18)$$

where $\dot{w}(t)$ can be zero before influx control (i.e., $a_{00}(t) = f(t)$), positive under control ($a_{00}(t) = \bar{a} < f(t)$), and negative after control (e.g., $a_{00}(t) = f(t) + \frac{w(t)}{\Delta t}$, where $\frac{w(t)}{\Delta t} > 0$ indicates the releasing rate of the unmatched requests; and $\Delta t$ is a pre-specified time interval).

Based on the above metrics, we now can measure the overall system performance by the total trip time spent in the service, $\mathcal{Z}$ (h), and the average, $\bar{\mathcal{Z}}$ (h/trip), as below:

$$\mathcal{Z} = \int_0^T w(t)dt + \int_0^T n_{01}(t)dt + \int_0^T n_{10}(t)dt \quad (19a)$$

$$\bar{\mathcal{Z}} = \frac{\mathcal{Z}}{\int_0^T f(t)dt} \quad (19b)$$

where the first term at RHS of (19a) is total waiting time spent by requests before being assigned; the second term means the waiting time for being picked up; and the third term is the total in-vehicle time.

## 2.3 Optimal control

To control the traffic in avoidance of gridlock, we proposed a density-based (DB) rule as follows:

*The DB rule*: At any time $t$, set $a_{00}(t) = \min\left(f(t) + \frac{w(t)}{\Delta t}, \bar{a}\right)$, if $\rho(t) \geq \rho_k$; else set $a_{00}(t) = \min\left(f(t) + \frac{w(t)}{\Delta t}, \frac{\rho_k - \rho(t)}{\Delta t} L\right)$, where $\rho_k$ is the critical density that generates the maximum flow, i.e., $\rho_k = \text{argmax}_{\rho(t)}\left(\rho(t)v(t, \rho(t))\right)$; and $\bar{a}$ is the control variable derived as below.

It is desired to keep the density not exceeding the critical density; otherwise, more delays would be incurred on all vehicles in the congested traffic. Thus, we have:

$$N(t) = n_{00}(t) + n_{01}(t) + n_{10}(t) \leq \rho_k L \tag{20}$$

Taking derivative with respect to $t$ on both sizes yields,

$$\dot{n}_{00}(t) + [a_{00}(t) - p_{01}(t)] + [p_{01}(t) - d_{10}(t)] \leq 0 \tag{21}$$

where $\dot{n}_{00}(t) = 0$ under perfect matching between supply and controlled demand. Therefore, (21) gives us

$$a_{00}(t) \leq \bar{a} \equiv d_{10}(t^*) = k_{10}(t^*, 0)v(t^*) \tag{22}$$

where $t^*$ is the time when $N^*(t^*) = \rho_k L$ or $\rho^*(t^*) = \rho_k$. The physical meaning of (22) is straightforward: The in-flux of vehicles should not exceed the maximum out-flux of the system.

The above DB rule normally functions as follows. For a typical peak period, $\rho(t)$ first rises with increasing $N(t)$. Before reaching $\rho_k$, $a_{00}(t)$ is not controlled and set to be $f(t)$ ($w(t) = 0$). Once $\rho(t)$ becomes no more less than $\rho_k$, $a_{00}(t)$ should be controlled by $\bar{a}$ (20) such that the system is kept stable at the maximum flow point; and the unmatched requests start accumulating ($w(t) \geq 0$). Afterward, when $f(t)$ falls down lower than $\bar{a}$, there is a room for the waiting requests entering the system, and the entering rate is controlled by $\frac{w(t)}{\Delta t}$. But overall, the total entering rate, i.e., $f(t) + \frac{w(t)}{\Delta t}$, should not exceed the accommodating capacity of the spared space in the network, which is measured by $\frac{\rho_k - \rho(t)}{\Delta t} L$.

The optimality of the above DB rule can be proved in a similar way of Daganzo (2007). For the completeness, we phrase the optimality theorem and proof as follows. Let $F(t) \equiv \int_0^t f(s)ds$ be the cumulative arrivals of travel requests, $A(t) \equiv \int_0^t a_{00}(s)ds$ be the cumulative assignment of requests, and $D(t) \equiv \int_0^t d_{10}(s)ds$ the cumulative departures of trips.

**Optimality Theorem**: For a network with a concave and unimodal speed-density relation, the DB rule determines a pair of curves $(A^*(t), D^*(t))$ that yields the minimum total system cost $(\mathcal{Z})$ given a demand pattern $F(t)$.

**Proof**: The total passenger time, i.e., the enclosed area between $F(t)$ and $D(t)$, cannot be further squeezed due to the fact that: (i) Under the state of $\rho(t) > \rho_k$, the $(A(t), D(t))$ is dictated by the exiting rate $\bar{a}(\rho_k)$ (20), which by definition generates the maximum exiting flow. For the state of $\rho(t) < \rho_k$, the $(A(t), D(t))$ is dictated by the possible in-flux (and also the left accommodating capacity), which cannot be enlarged either. The proof is completed. □

Note that the above DB rule bears similarity with the accumulation-based (AB) rule proposed by Daganzo (2007): Both rely on only monitoring the density to trigger the control, and thus are robust without the need of forecasting. Traditional technologies, e.g., magnetic loops at certain points of the network, work effectively for the two control rules. Furthermore, our DB rule indicates a more precise control variable $\bar{a}$ as opposed to that of Daganzo (2007), i.e., the endogenous input flow $dO$, of which the observation is not easy.

## 3. Solution method

From section 2, we have one ODE (1) and two PDEs (2-3) with three unknown variables, $n_{00}(t), k_{01}(t,x), k_{10}(t,x)$. The solutions can be found using finite difference method (FDM) as follows.

Let $T$ (h) and $X$ (km) be the study period and the maximum travel distance, respectively; discretize them into $J$ steps with $\Delta t = \frac{T}{J}$ and $I$ steps with $\Delta x = \frac{X}{I}$. The detailed steps of FDM is presented as below, where variables with subscripts or superscripts of $i, j$ denote the corresponding values of the $i$th distance, $i\Delta x$, and $j$th time step, $j\Delta t$.

Initialize zero values to $\{n_{00}^{i,0}, k_{01}^{i,0}, k_{10}^{i,0}, K_{01}^{i,0}, K_{10}^{i,0}\}_{i=0,1,\ldots,I}$

Proceed forward in time $j = 0, 1, \ldots, J$

    Continue if $N^{j+1} = 0$; otherwise, skip this time step and go to $j + 1$.

    Compute $v^j$ by (7), a variable time step $\Delta t_j = \frac{\Delta x}{v^j}$, and the proceeded steps $m = \frac{\Delta t_j}{\Delta t}$.

    Determine $a_{00}^j$ by the DB rule; set $p_{01}^j = k_{01}^{0,j} v^j$.

    Proceed backward in distance $i = I, I-1, \ldots, 0$

Update: $k_{01}^{i,j+m} = k_{01}^{i+1,j} + a_{00}^j \tilde{\varphi}_{01}(t_j, x_i)\Delta t_j$, and $k_{10}^{i,j+m} = k_{10}^{i+1,j} + p_{01}^j \tilde{\varphi}_{10}(t_j, x_i)\Delta t_j$.*

Update: $K_{01}^{i,j+m} = K_{01}^{i+1,j} + a_{00}^j \tilde{\Phi}_{01}(t_j, x_i)\Delta t_j$, and $K_{10}^{i,j+m} = K_{10}^{i+1,j} + p_{01}^j \tilde{\Phi}_{10}(t_j, x_i)\Delta t_j$.**

Update variables at $l \in (j, m)$ time steps: $k_{01}^{i,j+l} = k_{01}^{i,j} + \frac{\left(k_{01}^{i,j+m} - k_{01}^{i,j}\right)}{m-l}$, $k_{10}^{i,j+l} = k_{10}^{i,j} + \frac{\left(k_{10}^{i,j+m} - k_{10}^{i,j}\right)}{m-l}$ and $K_{01}^{i,j+l} = K_{01}^{i,j} + \frac{\left(K_{01}^{i,j+m} - K_{01}^{i,j}\right)}{m-l}$, $K_{10}^{i,j+l} = K_{10}^{i,j} + \frac{\left(K_{10}^{i,j+m} - K_{10}^{i,j}\right)}{m-l}$.†

End backward process in trip distance.

Update: $\dot{w}^{j+1} = f^{j+1} - a_{00}^{j+1}$; $w^{j+1} = \Delta t \sum_{l=0}^{j+1} \dot{w}^l$; $n_{00}^{j+1} = \min\{f^{j+1}\Delta t + w^j, \rho_k L - K_{01}^{0,j+1} - K_{10}^{0,j+1}\}$; and $N^{j+1}$.‡

End time process if $j$ reaches $J$.

* See the discrete models (5,6); ** They are the discrete versions of (10,11) obtained by integrating (5,6). † The equations yield smooth transitions of variables between time steps $j, m$. ‡ In this step, $\tilde{\varphi}_{01}(t_j, x)$ and $\tilde{\varphi}_{10}(t_j, x)$ may also be updated according to certain rules (See the section of numerical examples).

## 4. Extension to ride-pooling with participants $c > 1$

The above model and solution can be extended to ride-pooling services with the number of participants $c > 1$. In doing so, the operation strategy of ride-pooling needs to be specified. For instance, a particular operation scheme (as described in Daganzo, 1978; Daganzo and Ouyang, 2009) among others functions as follows.

- Ride requests are assembled by a central operator and assigned to the closest idle vehicles.
- Routing plan is made to vehicles that finished demand matching and takes no more new requests afterward.
- The active vehicle picks up all assigned passengers and then deliver them.

Under the above operation scheme, the numbers of three types of vehicles are revised to be $n_{00}(t) = \text{sum}\{n_{0j}(t)\}_{j<c}$, $n_{0c}(t) = \text{sum}\{n_{ij}(t)\}_{i+j=c, j>0}$, and $n_{c0}(t) = \text{sum}\{n_{i0}(t)\}_{i\leq c}$.[3] And for the active vehicles, the PDF functions of the desired travel distance

---

[3] Other operation schemes may have three types of vehicles overlapped when collecting vehicles accept new request assignment, or vehicles alternate pick-ups and deliveries (Daganzo, 1978). In these cases, the dynamic

distributions are, $\tilde{\varphi}_{0c}(t,x)$ and $\tilde{\varphi}_{c0}(t,x)$ (km$^{-1}$), respectively.[4] With these new variables, the above derivations for the dynamic models (9a-c) in section 2 still hold with an update to $a_{00}(t) = \frac{f(t)}{c}$ if no control, or $a_{00}(t) = \min\left(\frac{f(t)}{c}, \bar{a}\right)$ under a control.

To describe the dynamics of trips in the system, we define $h_{0c}(t,x)$ and $h_{c0}(t,x)$ as the densities of trips of active vehicles in two states with remaining travel distance being $x$ at $t$. Then, the following trips-oriented bathtub models can be formulated and solved together with the vehicles-oriented bathtub models:

$$\frac{\partial}{\partial t} h_{0c}(t,x) - \frac{\partial}{\partial x} h_{0c}(t,x) v(t) = c \cdot a_{00}(t) \tilde{\varphi}_{0c}(t,x) \tag{23}$$

$$\frac{\partial}{\partial t} h_{c0}(t,x) - \frac{\partial}{\partial x} h_{c0}(t,x) v(t) = c \cdot p_{0c}(t) \tilde{\varphi}_{c0}(t,x) \tag{24}$$

It is also noticed that the ride-pooling services allow us to eliminate the waiting time of unmatched requests by varying $c(t)$ such that $a_{00}(t) = \frac{f(t)}{c(t)} = \bar{a}$ if $\rho(t) \geq \rho_k$ and $\dot{w}(t) = f(t) - c(t) a_{00}(t) = 0$. Therefore, the optimal ride-pooling size $c(t)$ for the saturated states of $\rho(t) \geq \rho_k$ can be expressed by:

$$c(t) = \max\left(\left\lceil \frac{f(t)}{[k_{c0}(t^*,0) v(t^*)]_{t^*|\rho(t^*)=\rho_k}} \right\rceil, 1\right) \tag{25}$$

where $\lceil \cdot \rceil$ returns the ceiling integer of the argument.

The above setting of $c(t)$, however, may not the optimal solution that minimizes the total passenger time over the entire study period. This is because a tradeoff exists for the unsaturated state of $\rho(t) < \rho_k$: On one hand, an increase in $c(t)$ reduces the in-flux of vehicles in need, which contributes to a less congested traffic transporting the same number of passengers; on the other hand, the increase in $c(t)$ brings not only longer expected travel distance of active vehicles, which contributes to a more congested traffic, but also additional waiting time for matching the requests in a pool. Therefore, the following optimization model is formulated to find the optimal $c(t)$ under the unsaturated traffic condition.

$$\min_{c(t)} Z = \int_0^T \frac{c(t)}{2f(t)} f(t) dt + \int_0^T \frac{n_{0c}(t) c(t)}{2} dt + \int_0^T \frac{n_{c0}(t) c(t)}{2} dt \tag{26a}$$

subject to:

---

models may be formulated in a similar manner but with more complicated forms for vehicles of each particular state, $n_{ij}(t)$.

[4] They again can be calibrated using real data. However, these data of a particular $c$-rides-pooling service may not be rich if available. So, we reckon that $\tilde{\varphi}_{0c}(t,x) = \frac{\tilde{\varphi}_{01}(t,x)}{c}$ and $\tilde{\varphi}_{c0}(t,x) = \frac{\tilde{\varphi}_{10}(t,x)}{\sqrt{c}}$ holds approximately. The underlining logic is that given demand-supply conditions, the expected collecting and delivering tour distances of a $c$-rides-pooling vehicle are $c$ and $\sqrt{c}$ times longer than those of a single-ride vehicle (See Daganzo, 1984), which stretches the distributions and correspondingly shrinks PDFs $c$ and $\sqrt{c}$ times.

$$c(t) \in \{1,2,\dots,c^{\max}\}, \forall t \tag{26b}$$

where the first term in RHS of (26a) accounts for the time for matching $c(t)$ ride-pooling participants, and $\frac{c(t)}{2f(t)}$ is the average waiting time per request during the matching; and $c^{\max}$ is the maximum ride-pooling size (e.g., 4 for a sedan, 7 for a minivan, and 15 for a minibus).

The problem (26) can be written in a discrete form:

$$\min_{c^j} \mathcal{Z} = \sum_{j=0}^{J} \frac{\left(1+n_{0c}^j+n_{c0}^j\right)c^j}{2}\Delta t \tag{27a}$$

subject to:

$$c^j \in \{1,2,\dots,c^{\max}\}, \forall j \tag{27b}$$

The above problem (27) is a sequential decision problem and thus can be solved by dynamic programming (DP). Main steps are briefed as follows. Let $S^j \equiv \{n_{00}^j, n_{0c}^j, n_{c0}^j\}$ be the system state at $t^j$, the original problem (27) can be rewritten into an equivalent DP problem:

$$\mathcal{Z}(S^j) = \min_{c^j} z^j(S^j) + \mathcal{Z}\left(S^{j+1}(S^j,c^j,f^j)\right) \tag{28a}$$

subject to:

$$S^{j+1} = S^{j+1}(S^j,c^j,f^j) \tag{28b}$$

$$\mathcal{Z}(S^{J+1}) = 0 \tag{28c}$$

$$c^j \in \{1,2,\dots,c^{\max}\}, \forall j \tag{28d}$$

where $z^j = \frac{\left(1+n_{0c}^j+n_{c0}^j\right)c^j}{2}\Delta t$ is the system cost of onetime step at $t^j$; (28b) is the state transition function; and (28c) is the boundary condition.

The DP problem (28) can be solved by the classic backward recursion for $j = J, J-1, \dots, 0$, and $\mathcal{Z}(S^0)$ is the optimal objective function value of (27a). Then, the optimal $c^{j*}$ can be traced in a forward way for $j = 0,1,2,\dots,J$ in correspondence to the state $S^{j*} = S^j(S^{j-1}, c^{j-1*}, f^{j-1})$ that have the lowest $\mathcal{Z}(S^{j*})$:

$$c^{j*} = \mathrm{argmin}_{c^j}\, z^j(S^{j*}) + \mathcal{Z}\left(S^{j+1}(S^{j*},c^j,f^j)\right) \tag{29}$$

Specifically, considering the possible states of the system are tremendous but should obey to the traffic dynamics models, the Monte Carlo method can be used to generate trajectories of $\{S^j\}_{j=0,1,\dots,J}$ under the range of $c^j$. The detailed algorithm steps can be found in any DP textbooks (e.g., Bellman, 1957; Bertsekas, 2007) and thus omitted here.

## 5. Numerical examples

### 5.1 Basic settings

The following settings are adapted from the example in Jin (2020). Consider a small town with area of $\mathcal{A} = 5$ km². Her network is with $L = 10$ lane-kms. The speed-density relation is expressed by

$$V(\rho) = \min\left\{30, \frac{750}{\rho}, 10\left(\frac{200}{\rho} - 1\right)\right\},$$

of which the critical density is known to be $\rho_k = 125$ vehicles/km. The demand during a peak period is given as

$$f(t) = \max\{0, \min\{100t, 100, 100(1-t)\}\} \text{ (trips/h)}.$$

The desired travel distances of ride-sourcing vehicles are assumed to follow uniform distributions[5] as below:

$$\tilde{\varphi}_{0c}(t,x) = \max\left(0, 1 - \frac{x}{2\tilde{B}_{0c}(t)}\right), c = 1,2,\ldots$$

$$\tilde{\varphi}_{c0}(t,x) = \max\left(0, 1 - \frac{x}{2\tilde{B}_{c0}(t)}\right), c = 1,2,\ldots$$

where $\tilde{B}_{01}(t)$ and $\tilde{B}_{10}(t)$ are the expected travel distances of collecting and delivering vehicles. According to Daganzo (1984), under the assumption of the spatially uniform distribution of trips' origins and destinations, they can be expressed by:

$$\tilde{B}_{01}(t) = \hbar\sqrt{\frac{\mathcal{A}}{n_{00}(t)}}, \tilde{B}_{0c}(t) = c\hbar\sqrt{\frac{\mathcal{A}}{n_{00}(t)}}$$

$$\tilde{B}_{10}(t) = \hbar'\sqrt{\mathcal{A}}, \tilde{B}_{c0}(t) = \hbar'\sqrt{\mathcal{A}c}$$

where in Manhattan metrics $\hbar$ and $\hbar'$ take 0.63 and 1.15, respectively.

We next demonstrate the traffic dynamics of non-shared taxis and ride-pooling vehicles with and without control, respectively.

### 5.2 Dynamics of non-shared taxis

Figure 1a-d illustrates the traffic dynamics of uncontrolled non-shared taxis with respect to the speed, characteristic travel distance, and numbers of active vehicles with ranges of remaining travel distances. It is observed in Figure 1a, b that the traffic reaches gridlock (i.e., $v(t) = \dot{z}(t) = 0$) soon at time 37.4 mins. Figure 1c shows that the travel distance distribution of collecting vehicles shrinks as time goes because of the increasing number of available idle

---

[5] Other distributions may also apply. For instance, in studies on traveling salesman problems (TSP) or dial-a-ride problems (DAR), some found that the TSP tour length follows a normal distribution (Vinel and Silva, 2018).

vehicles. From Figure 1d, it is seen the cumulation of vehicles intensifies at around 28 mins and quickly falls into the gridlock state with the majority of shorter remaining distances.

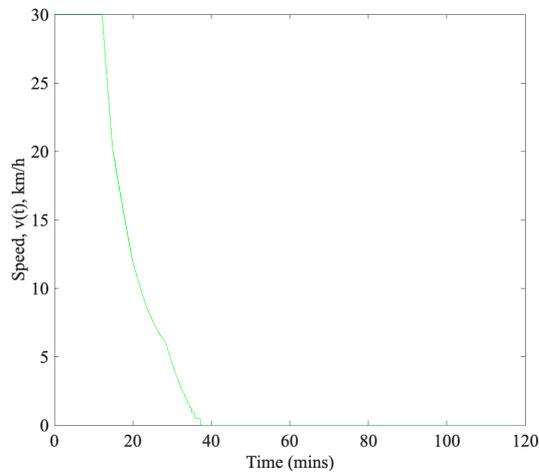

a. Speed, $v(t)$ (km/h)

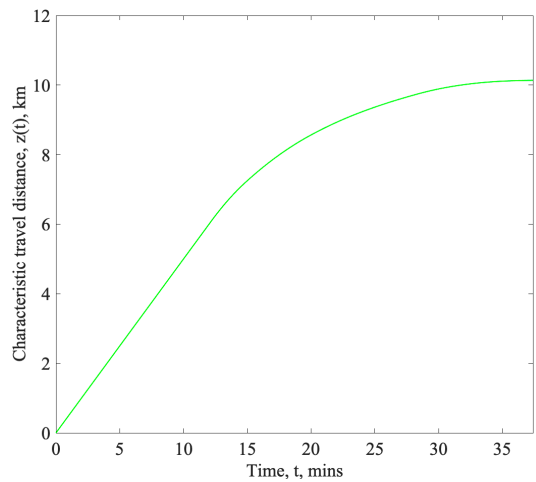

b. Characteristic travel distance, $z(t)$ (km)

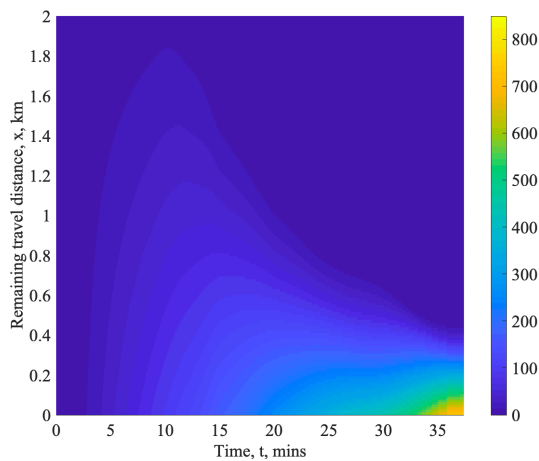

c. Collecting vehicles, $K_{01}(t,x)$ (vehicles)

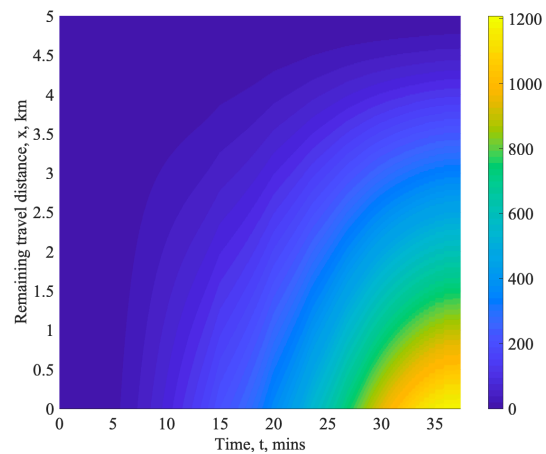

d. Delivering vehicles, $K_{10}(t,x)$ (vehicles)

**Figure 1. Dynamics of non-shared taxis without control**

For comparison, Figure 2a-d depicts the dynamics of controlled non-shared taxis. Evidently, the gridlock is avoided. The speed in Figure 2a is stabilized at around 6 km/h during the time interval [28, 82] mins with a duration about 54 mins. The control is also reflected in Figure 2c, where a narrow gap (about 2-mins long), indicating a short shutdown of in-flow, can be observed before the time of 40 mins. At that time, the accumulation of delivering vehicles is critically intensified; see the interval of [30, 40] mins in Figure 2d. After the implementation of the control, both distributions of $K_{10}(t,x)$ and $K_{01}(t,x)$ are stabilized until the incoming demand starts fading. Under the control, all trips are accomplished and the average time per trip is 37.18 mins.

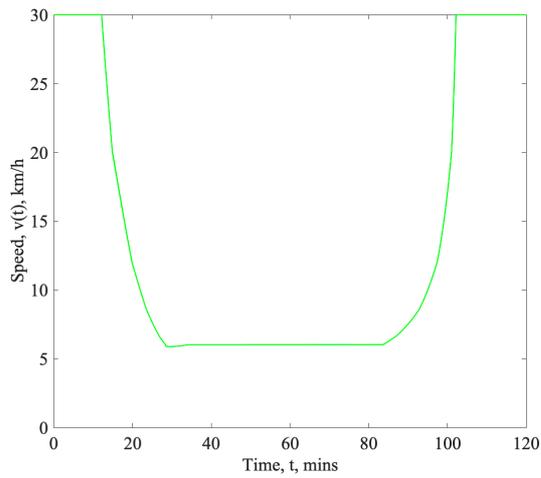

a. Speed, $v(t)$ (km/h)

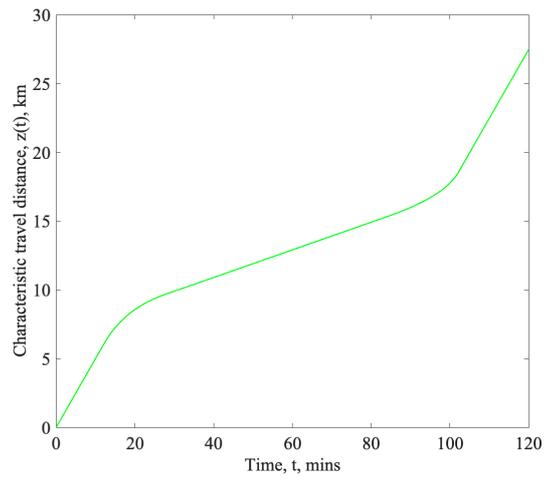

b. Characteristic travel distance, $z(t)$ (km)

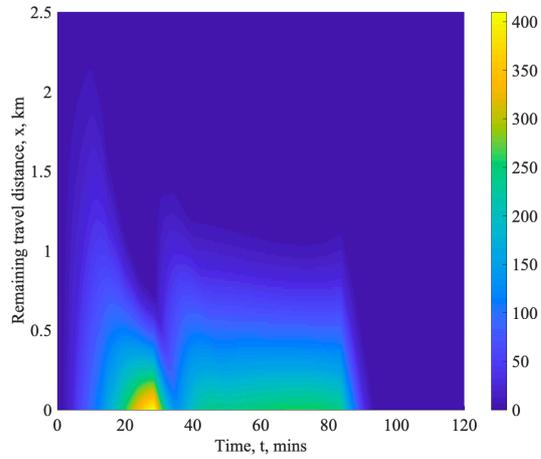

c. Collecting vehicles, $K_{01}(t,x)$ (vehicles)

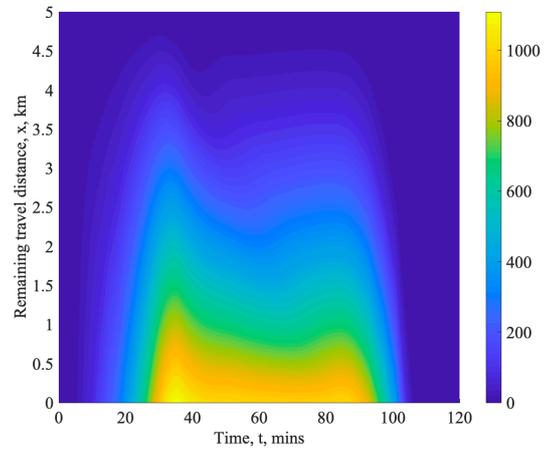

d. Delivering vehicles, $K_{10}(t,x)$ (vehicles)

**Figure 2. Dynamics of non-shared taxis with the optimal control**

The control inevitable results some requests having to wait before being assigned to ride-sourcing cars. Figure 3a shows the changes of total waiting requests, which reflects the control on the inflow. The information can also be represented using the cumulative curves of trip requests and those being served in active vehicles, as in Figure 3b. From Figure 3b, it can be told that the vertical distance between the green and blue curves indicates the instant number of waiting requests (including those waiting to be assigned and those waiting to be picked); and the horizontal distance means the waiting time for the request entering the system at $t$. Thus, the enclosed area denotes the total waiting time. Similarly, the instant number of trips being delivered, the travel time for trips entering at certain time, and the total trip times can also be conveyed from Figure 3b.

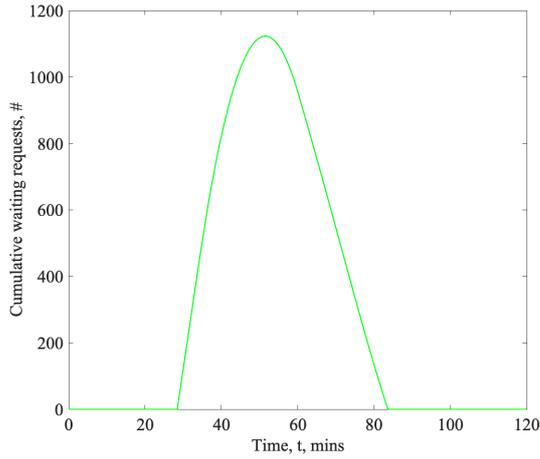
a. Total waiting requests, $w(t)$

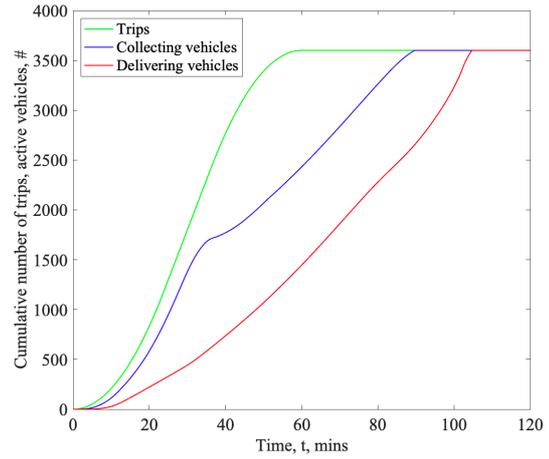
b. Cumulative trips and active vehicles

**Figure 3. Total waiting requests and cumulative trip requests, active vehicles**

### 5.3 Dynamics of ride-pooling vehicles

Several scenarios are experimented for the controlled ride-pooling services with ride-pooling sizes capped by $c(t) \leq 2$, $c(t) \leq 3$, $c(t) \leq \inf$, and dictated by the optimal $c^{j*}$ from (29), respectively. Figure 4a-d presents the distributions of active vehicles. Comparatively, the critical peak period of traffic (see the highlighted area in Figure 4a-c) is shortened by larger-sized ride-pooling vehicles, although they are more scattered in travel distances. Further observation finds that in the first three scenario the ride-pooling size rises from 1 to the maximum at the same critical time (around 28 mins); however, in order to prevent gridlock happening, the in-flow control in scenarios with $c(t) \leq 2,3$ has to also function and makes unassigned requests waiting; and it is not necessary in scenario 3, of which the maximum ride-pooling size is as high as 7 passengers/vehicle. Lastly, Figure 4d shows the dynamics under the optimal $c^{j*}$, which is found to be 2 passengers/vehicle. As seen, the critical peak period is the shortest and no requests waiting for assignment. The average time per trip is reduced to 31.21 mins/trip.

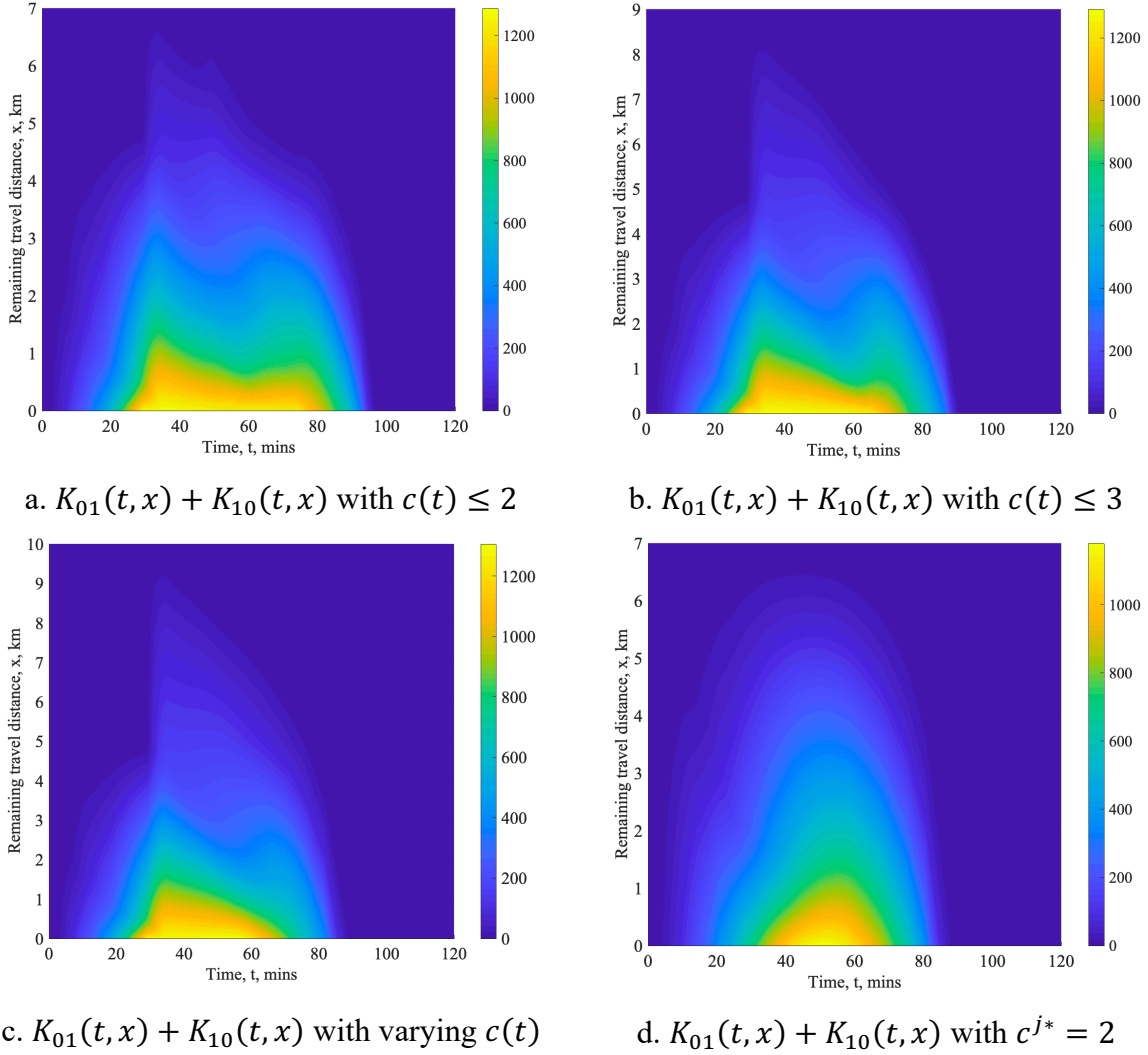

a. $K_{01}(t,x) + K_{10}(t,x)$ with $c(t) \leq 2$

b. $K_{01}(t,x) + K_{10}(t,x)$ with $c(t) \leq 3$

c. $K_{01}(t,x) + K_{10}(t,x)$ with varying $c(t)$

d. $K_{01}(t,x) + K_{10}(t,x)$ with $c^{j*} = 2$

**Figure 4. Effects of ride-pooling sizes**

## 6. Conclusions

This paper proposes an interactive bathtub model for describing the dynamics of ride-sourcing vehicles in the process of demand-vehicle matching, picking up, and delivering. The model is parsimonious and requests only information of the total lane length of the network, the in-flux of demand, and the travel distance distributions of ride-sourcing vehicles, which are all obtainable with low cost from a practical perspective. The proposed model captures the dynamic states of not only active vehicles in traffic but also trips in the services. To manage the traffic efficiently and avoid gridlock, a density-based (DB) rule is proposed and proved to be the optimal control strategy. Precise control variable is derived.

The above model and control rule also apply to ride-pooling services under certain operation scheme. If the ride-pooling size is allowed to vary, we derive the optimal size under the saturated traffic state; whereas for the unsaturated states, the optimal ride-pooling size needs to be found by solving an optimization problem. For the latter, we formulate into a dynamic programming problem and accordingly provide the solution algorithm. Numerical examples demonstrate the effectiveness of the control and the varying ride-pooling sizes on managing the traffic.

The proposed model can be used by the authorities to monitor, predict, and accordingly control (via TNCs) the traffic; or by the TNCs to determine the fleet size in advance, dynamically dispatch vehicles, and macroscopically measure the productivity of the service.

Admitted, the proposed model is much simplified as compared to the real world, where trips by private cars and ride-sourcing services co-exist, and the ride-sourcing services may also be mixed by non-shared taxis and ride-pooling cars with various sizes. The perfect instant matching between requests and vehicles may also not hold in reality; and the supply of ride-sourcing vehicles may be lagged behind real-time requests and also capped by a maximum fleet size. Additionally, in response to the control, travelers may alter their behaviors by quitting the queue of unmatched requests upon a maximum waiting time or directly taking the public transit. The TNCs may thus be reluctant to obey the control rule under the burden of competing for customers. Even so, common interests indeed exist among all stakeholders (i.e., ordinary travelers, TNCs, and the authorities) on a healthy traffic state, which indicates low delays, high productivity, and large social welfare. Some of these issues may be addressed in the future research, e.g., by introducing pricing instrument, travelers' choice models, and developing multimodal bathtub models.

**Acknowledgments**

The research was supported by a fund provided by the National Natural Science Foundation of China (No. 51608455) and Sichuan Science & Technology Program (No. 2020YFH0038).

## Appendix A: Notations

**Table A1 List of notations**

| Variables | Definitions | Units |
|---|---|---|
| $A(t)$ | Cumulative number of assigned/matched requests at $t$ | # |
| $D(t)$ | Cumulative number of departed/delivered requests at $t$ | # |
| $K_{01}(t,x), K_{10}(t,x)$ | Numbers of collecting, delivering non-shared taxis at $t$ with remaining travel distances no less than $x$ | # |
| $F(t)$ | Cumulative number of requests at $t$ | lane-km |
| $I, J, i, j, \Delta x, \Delta t$ | Space-time discretization variables | |
| $L$ | Total lane length of the network | |
| $N(t)$ | Total number of active vehicles in the network at $t$ | # |
| $S^j$ | State variable of the system at $j$th time step | |
| $T$ | Time duration of the study period | h |
| $X$ | Maximum travel distance | km |
| $\mathcal{Z}, \bar{\mathcal{Z}}, z^j$ | The total, average, and $j$th time interval system cost | pax·h |
| $a_{00}(t)$ | Rate of vehicles changing from idle to collecting state | veh/h |
| $\bar{a}$ | Control variable | veh/h |
| $c(t)$ | Ride-pooling size, if fixed denoted as $c$ | pax/veh |
| $d_{10}(t), d_{c0}(t)$ | Rate of vehicles changing from delivering to idle state | veh/h |
| $f(t)$ | In-flux of requests for ride-sourcing services | trips/h |
| $h_{0c}(t,x), h_{c0}(t,x)$ | Densities of trips being collected and delivered at $t$ with remaining travel distance being $x$ | trips/km |
| $k_{01}(t,x), k_{10}(t,x)$ $k_{0c}(t,x), k_{c0}(t,x)$ | Densities of collecting and delivering vehicles at $t$ with remaining travel distance being $x$ | veh/km |
| $n_{00}(t), n_{01}(t), n_{10}(t)$ | The numbers of idling, collecting, and delivering non-shared taxis at $t$ | veh |
| $p_{01}(t), p_{0c}(t)$ | Rate of vehicles changing from collecting to delivering state | veh/h |
| $s(t)$ | Supply rate of ride-sourcing vehicles | veh/h |
| $v(t)$ | Travel speed in the network at $t$ | km/h |
| $w(t)$ | Total number of unmatched requests at $t$ | # |
| $\mathcal{A}$ | Area of the study domain | km² |
| $\tilde{B}_{01}(t), \tilde{B}_{10}(t)$ $\tilde{B}_{0c}(t), \tilde{B}_{c0}(t)$ | Expected travel distances of collecting and delivering vehicles at $t$ | km |
| $\mathcal{k}, \mathcal{k}'$ | Parameters in models of expected travel distances | |
| $\mathcal{n}_{00}(t), \mathcal{n}_{0c}(t), \mathcal{n}_{c0}(t)$ | The numbers of idling, collecting, and delivering ride-pooling vehicles at $t$ | # |
| $\rho(t), \rho_k$ | The average, critical density of vehicles in the network at $t$ | veh/km |

| $\tilde{\varphi}_{01}(t,x), \tilde{\varphi}_{10}(t,x)$ $\tilde{\varphi}_{0c}(t,x), \tilde{\varphi}_{c0}(t,x)$ | PDFs of collecting, delivering vehicles entering at $t$ with desired travel distance being $x$ | km$^{-1}$ |
|---|---|---|
| $\widetilde{\Phi}_{01}(t,x), \widetilde{\Phi}_{10}(t,x)$ | CDFs of collecting, delivering vehicles entering at $t$ with desired travel distances no less than $x$ | |
| $\varphi_{01}(t,x), \varphi_{10}(t,x)$ | PDFs of collecting, delivering vehicles at $t$ with remaining travel distance being $x$ | km$^{-1}$ |
| $\Phi_{01}(t,x), \Phi_{10}(t,x)$ | CDFs of collecting, delivering vehicles at $t$ with remaining travel distances no less than $x$ | |